# A BERRY–ESSEEN THEOREM FOR SAMPLE QUANTILES UNDER WEAK DEPENDENCE[1]

By S. N. Lahiri and S. Sun

*Texas A&M University and Wright State University*

This paper proves a Berry–Esseen theorem for sample quantiles of strongly-mixing random variables under a *polynomial* mixing rate. The rate of normal approximation is shown to be $O(n^{-1/2})$ as $n \to \infty$, where $n$ denotes the sample size. This result is in sharp contrast to the case of the sample mean of strongly-mixing random variables where the rate $O(n^{-1/2})$ is not known even under an exponential strong mixing rate. The main result of the paper has applications in finance and econometrics as financial time series data often are heavy-tailed and quantile based methods play an important role in various problems in finance, including hedging and risk management.

**1. Introduction.** Sample quantiles of time series data play an important role in robust statistical inference about various process parameters, particularly when the underlying distribution is heavy-tailed or when presence of outliers is suspected [cf. Dutter, Filzmoser, Gather and Rousseeuw (2003)]. Although asymptotic normality of the sample quantiles under dependence is known, accuracy of the corresponding normal approximation has remained largely unexplored. In this paper, we establish a Berry–Esseen theorem for the sample quantile with the optimal rate $O(n^{-1/2})$ for a large class of weakly dependent time series. Apart from its foundational role in statistical inference for time series data, the Berry–Esseen result of the paper also has important applications in finance and econometrics. It is well known [cf. Mittnik and Rachev (2001)] that financial time series data often are heavy-tailed. As a result, quantile based methods are being increasingly developed and employed in diverse problems in finance, such as, quantile-hedging [cf. Föllmer and Leukert (1999)], optimal portfolio allocation [cf. Dmitrašinović-Vidović and Ware (2006)], risk management [cf. Melnikov and Romaniuk

Received May 2007; revised March 2008.

[1]Supported by NSF Grants DMS-03-06574 and DMS-07-42690.

(2006)], and so forth. The recent surge in interest in using quantile based methods in finance and related areas calls for a better understanding of the theoretical properties of the sample quantiles of time series data in greater generality. The main result of the paper takes a step in this direction and establishes the optimal rate in the Berry–Esseen theorem for a large class of weakly dependent processes that require a *polynomial strong mixing* condition.

To describe the result of the paper, let $\{X_i\}_{i \in \mathbb{Z}}$ be a sequence of stationary random variables with strong mixing coefficient $\alpha_X(n) = \sup\{|P(A \cap B) - P(A)P(B)| : A \in \mathcal{F}_{-\infty}^k, B \in \mathcal{F}_{k+n}^\infty, k \in \mathbb{Z}\}$, where $\mathcal{F}_a^b = \sigma\langle X_i : i \in [a,b] \cap \mathbb{Z}\rangle$ is the $\sigma$-field generated by $\{X_i : i \in [a,b] \cap \mathbb{Z}\}$, $-\infty \leq a \leq b \leq \infty$ and where $\mathbb{Z} \equiv \{0, \pm 1, \pm 2, \ldots\}$ denotes the set of all integers. Let $F$ denote the distribution function (d.f.) of $X_1$, that is, $F(x) = P(X_1 \leq x), x \in \mathbb{R}$. For $p \in (0,1)$, let

$$F^{-1}(p) = \inf\{x : F(x) \geq p\} \tag{1.1}$$

denote the $p$th quantile of $F$. An estimator of the population parameter $F^{-1}(p)$ is given by the sample $p$th quantile

$$F_n^{-1}(p) = \inf\{x : F_n(x) \geq p\}, \tag{1.2}$$

where $F_n(x) = n^{-1} \sum_{i=1}^n I(X_i \leq x), x \in \mathbb{R}$, denotes the empirical distribution function (e.d.f.) of $X_1, \ldots, X_n$ and where $I$ denotes the indicator function, with $I(S) = 0$ or $1$ according as the statement $S$ is false or true. When the process $\{X_i\}_{i \in \mathbb{Z}}$ is strongly mixing at a (suitable) polynomial rate [i.e., $\alpha_X(n) = O(n^{-a})$ as $n \to \infty$ for some suitable $a \in (0, \infty)$] and $F$ is differentiable at $F^{-1}(p)$ with a positive derivative $f(F^{-1}(p)) > 0$, it is known [see, e.g., Sen (1972), Sun and Lahiri (2006)] that

$$\sqrt{n}(F_n^{-1}(p) - F^{-1}(p)) \to^d N\left(0, \frac{\sigma_\infty^2(F^{-1}(p))}{f^2(F^{-1}(p))}\right) \tag{1.3}$$

as $n \to \infty$, where $\sigma_\infty^2(x) = \sum_{i \in \mathbb{Z}} \text{Cov}(I(X_1 \leq x), I(X_{i+1} \leq x))$, $x \in \mathbb{R}$.

The main result of this paper refines (1.3) by specifying the rate of normal approximation to the distribution of $\sqrt{n}(F_n^{-1}(p) - F^{-1}(p))$. More precisely, it is shown that if the process $\{X_i\}_{i \in \mathbb{Z}}$ is strongly mixing at a certain *polynomial* rate and if the regularity conditions set forth in Section 2 hold, then

$$\sup_{x \in \mathbb{R}} |P(\sqrt{n}(F_n^{-1}(p) - F^{-1}(p)) \leq x) - \Phi(x/\tau_\infty(p))| = O(n^{-1/2}) \tag{1.4}$$

as $n \to \infty$,

where $\tau_\infty^2(p) = \sigma_\infty^2(F^{-1}(p))/f^2(F^{-1}(p))$, and where $\Phi$ denotes the d.f. of a standard normal variate. Thus, the Berry–Esseen theorem holds for the *sample quantile* of strongly mixing random variables under the conditions of Section 2. This is in marked contrast to the case of the *sample mean* of strongly mixing random variables, where a Berry–Esseen theorem with the rate $O(n^{-1/2})$ of normal approximation is not available. The best known rate for sums of strongly mixing random variables with an *exponentially* decaying mixing coefficient is only $O(n^{-1/2}(\log n)^c)$ for some suitable $c > 0$ [cf. Tikhomirov (1980), Dasgupta (1988)]. The rate $O(n^{-1/2})$ for the sample mean is available either when the $X$-process satisfies certain *stronger* forms of dependence conditions, like $\phi$-mixing [cf. Donkhan (1994)] or when the distance between the two probability distributions is measured using certain *smooth* metrics. Rio (1996) obtained the $O(n^{-1/2})$ rate for the sample mean under a suitable *uniform mixing* condition, which is known to be stronger than strong mixing [cf. Donkhan (1994)]. For certain smooth metrics, Utev (1991) established the same optimal rate for scaled sums of Banach space valued random elements under $\phi$-mixing. For strongly mixing random vectors, the $O(n^{-1/2})$ rate for the sample mean (under suitable smooth metrics) follow from the results of Götze and Hipp (1983) and Lahiri (1993) under an exponential mixing rate, and from those of Lahiri (1996) under a polynomial mixing rate. Recently, Bentkus and Sunklodas (2007) and Sunklodas (2007) have established the $O(n^{-1/2})$ rate for the sample mean under different smooth metrics, for both strongly mixing random variables and strongly mixing random fields. However, for most statistical applications, approximations to the distribution function is needed and the best known rate for approximation to the distribution function of the sample mean for strongly mixing random variables is still $O(n^{-1/2}(\log n)^c)$ for some $c > 0$.

Although the validity of the Berry–Esseen theorem with rate $O(n^{-1/2})$ for the sample mean of strongly mixing random variables remains unsolved, the main result of this paper establishes the desired optimal rate $O(n^{-1/2})$ for the sample quantiles in the strong mixing case, requiring only a polynomial decay of the mixing coefficient. In particular, the Berry–Esseen theorem of the paper extends the results of Reiss (1974) who establishes the $O(n^{-1/2})$ rate of normal approximation to the distributions of sample quantiles under independence. The proof of the main result here makes use of some arguments developed by Götze and Hipp (1983) and Lahiri (1993, 1996) for deriving Edgeworth expansions for sums of strongly mixing random variables and also crucially exploits properties of the probability integral transform $F_n^{-1}(\cdot)$ of the e.d.f. $F_n$. We also consider some important classes of dependent processes and verify the regularity conditions used in the paper.

The rest of the paper is organized as follows. In Section 2, we state the conditions, verify these for different classes of weakly dependent processes and state the main result. The proof of the main result along with some auxiliary lemmas are given in Section 3.

## 2. Main results.

2.1. *Conditions.* We prove the Berry–Esseen bound under a general framework introduced by Götze and Hipp (1983) in their seminal paper on asymptotic expansions for sums of weakly dependent random vectors. Suppose that the random variables $\{X_i : i \in \mathbb{Z}\}$ are defined on a probability space $(\Omega, \mathcal{F}, P)$ and that $\{\mathcal{D}_i : i \in \mathbb{Z}\}$ is a collection of sub-$\sigma$-fields of $\mathcal{F}$. For $-\infty \leq a \leq b \leq \infty$, let $\mathcal{D}_a^b = \sigma\langle\{\mathcal{D}_i : i \in [a,b] \cap \mathbb{Z}\}\rangle$ denote the smallest $\sigma$-field containing $\{\mathcal{D}_i, a \leq i \leq b, i \in \mathbb{Z}\}$. Also, let $\xi_p = F^{-1}(p)$. Recall that $\sigma_\infty^2(x) \equiv \sum_{i \in \mathbb{Z}} \text{Cov}(I(X_1 \leq x), I(X_{i+1} \leq x))$, $x \in \mathbb{R}$. Let $\mathbb{N} = \{1, 2, \ldots\}$ denote the set of all positive integers. We shall make use of the following conditions:

(C.1)     (i) $F$ is differentiable at $\xi_p$ with derivative $f(\xi_p) \in (0, \infty)$.
         (ii) $\sigma_\infty^2(\xi_p) \in (0, \infty)$.

(C.2) There exist constants $d \in (0,1)$ and $\alpha_0 > 12$ such that for all $n \geq 1$,

$$\alpha(n) \equiv \sup\{|P(A \cap B) - P(A)P(B)| : A \in \mathcal{D}_{-\infty}^i, B \in \mathcal{D}_{i+n}^\infty, i \in \mathbb{Z}\} \tag{2.1}$$
$$\leq d^{-1} n^{-\alpha_0}.$$

(C.3) There exist constants $d \in (0,1)$ and $\beta_0 > 12$ and random variables $X_{i,n}^\dagger$, $i \in \mathbb{Z}$, $n \in \mathbb{N}$ such that $X_{i,n}^\dagger$ is $\mathcal{D}_{i-n}^{i+n}$-measurable and

$$\beta(n) \equiv E|X_i - X_{i,n}^\dagger| \leq d^{-1} n^{-2\beta_0} \qquad \text{for all } i \in \mathbb{Z}, n \in \mathbb{N}. \tag{2.2}$$

(C.4) There exist constants $d \in (0,1)$ and $\gamma_0 > 12$ such that for all $m, n, r \in \mathbb{N}$ and $A \in \mathcal{D}_{r-m}^{r+m}$,

$$|P(A|\mathcal{D}_j : j \neq r) - P(A|X_j : 0 < |r - j| \leq m + n)| \leq d^{-1} n^{-\gamma_0}. \tag{2.3}$$

(C.5) There exist a constant $d \in (0,1)$ and sub-$\sigma$-fields $\mathcal{C}_i, i \in \mathbb{Z}$, of $\mathcal{F}$ such that for every $i \in \mathbb{Z}$, $\sigma\langle \mathcal{D}_j : j \neq i\rangle \cup \sigma\langle\{X_j : j \neq i\}\rangle \subset \mathcal{C}_i$ and

$$P(G_i(\xi_p) = 1) \leq p - d \tag{2.4}$$

where $G_i(y) = P(X_i \leq y | \mathcal{C}_i)$, $y \in \mathbb{R}$.

We now comment on the conditions. Condition (C.1) is a standard condition that is frequently used to ensure a nondegenerate limit distribution of the $p$th sample quantile under dependence. In the independent case, (C.1)(i) is also known to be necessary; see Lahiri (1992). Conditions (C.2)–(C.4) are similar to the conditions introduced in Götze and Hipp (1983) for deriving asymptotic expansion for sums of weakly dependent random vectors, where the right-hand sides of (2.1)–(2.3) were assumed to be *exponentially* decaying functions of $n$. The reduction to the polynomial rate here heavily relies on Lahiri (1996) which extends Götze and Hipp's (1983) results allowing polynomial decay of the coefficients in (2.1)–(2.3). Condition (C.2) is a strong mixing condition on the auxiliary $\sigma$-fields $\mathcal{D}_j$'s, which together with condition (C.3), imposes an approximate strong-mixing structure to the given random variables $\{X_i\}_{i \in \mathbb{Z}}$. Condition (C.4) is an approximate Markov condition and in particular, it is satisfied if $\{X_i\}_{i \in \mathbb{Z}}$ is an $m$th order Markov process for a fixed $m \in \mathbb{N}$. Condition (C.5) is a key regularity condition that perhaps needs some detailed discussion. To get some insight into condition (C.5), first suppose that the $X_i$'s are independent. In this case, if we take $\mathcal{D}_j = \sigma\langle X_j \rangle, j \in \mathbb{Z}$, and $\mathcal{C}_j = \sigma\langle \{X_i : i \neq j\} \rangle$, then $G_i(\xi_p) = P(X_i \leq \xi_p | \mathcal{C}_i) = P(X_1 \leq \xi_p) = p$, making the probability on the left-hand side of (2.4) zero, and hence, condition (C.5) easily holds. For the dependent case, however, $G_i(\xi_p)$ is a random quantity. In this case, if $P(G_i(\xi_p) = 1) = p$, then one can show that the conditional distribution of $I(X_i \leq \xi_p)$ given $\mathcal{C}_i$ is degenerate at 1 on a set of probability $p$ while it is degenerate at 0 on the complementary set of probability $q = 1 - p$. As a result, the conditional characteristic function of $I(X_i \leq \xi_p)$ given $\mathcal{C}_i$ becomes identically equal to 1 in absolute value on all of $\Omega$, and our bound on the factorized conditional characteristic function of the scaled sum $n^{-1/2} \sum_{j=1}^n (I(X_i \leq \xi_p) - F(\xi_p))$ no longer provides a useful estimate for the discrepancy between the d.f.s of the sample quantile and the limiting normal distribution. However, once this degeneracy is ruled out by condition (C.5), it is possible to derive a suitably small upper bound on the conditional characteristic function of the scaled sums $n^{-1/2} \sum_{j=1}^n (I(X_i \leq y) - F(y))$ uniformly over $y$ in a neighborhood of $\xi_p$ (cf. Lemma 3.3 in Section 3). We exploit some basic properties of the quantile function in conjunction with this bound to establish the $O(n^{-1/2})$-order bound for the sample quantiles.

As in Götze and Hipp (1983), the $\sigma$-fields $\mathcal{D}_j$'s and $\mathcal{C}_j$'s are introduced to add flexibility in the verification of conditions (C.2)–(C.5). Below we consider some important examples and choose the $\sigma$-fields $\mathcal{D}_j$'s and $\mathcal{C}_j$'s suitably to show that condition (C.5) is quite unrestrictive.

2.2. *Examples.*

EXAMPLE 2.1. Suppose $\{X_i\}_{i \in \mathbb{Z}}$ is $m$-dependent for some $m \in \mathbb{Z} \cup \{0\}$, that is, $\sigma\langle \{X_i : i \leq k\} \rangle$ and $\sigma\langle \{X_i : i \geq k + m + 1\} \rangle$ are independent for all $k \in \mathbb{Z}$. Then we take $\mathcal{D}_j = \sigma\langle X_j \rangle$ and $\mathcal{C}_j = \sigma\langle X_i : i \neq j \rangle$, $j \in \mathbb{Z}$. Then it is easy to check that conditions (C.2)–(C.4) hold with $X_{i,m}^{\dagger} = X_i$ for all $i \in \mathbb{Z}$, $m \in \mathbb{N}$ and with $\alpha_0, \beta_0, \gamma_0$ arbitrarily large. Furthermore, in this case, condition (C.5) reduces to

(2.5) $$P(P(X_0 \leq \xi_p | X_i : 0 < |i| \leq m) = 1) < p.$$

Let $G_0$ denote the conditional distribution of $X_0$ given $\{X_i : 0 < |i| \leq m\}$. Suppose that there exist a set $A \in \mathcal{F}$ with $P(A) > 0$ and real numbers $\varepsilon, a, b$ with $\varepsilon \in (0, 1/2)$, $a \leq \xi_p \leq b$ such that $G_0$ puts at least $\varepsilon$ mass on $(a - \varepsilon, a]$ and on $(b, b + \varepsilon]$ on the set $A$, that is, if

(2.6) $\quad G_0((a - \varepsilon, a]) > \varepsilon, \qquad G_0((b, b + \varepsilon]) > \varepsilon \qquad$ for all $\omega \in A$.

We claim that (2.5) holds if (2.6) holds. To see this, note that (writing $G_0$ also to denote the distribution function),

(2.7)
$$p = F(\xi_p) = EG_0(\xi_p)$$
$$= P(G_0(\xi_p) = 1) + EG_0(\xi_p)I(0 < G(\xi_p) < 1)$$

so that $P(G_0(\xi_p) = 1) \leq p$. If possible, now suppose that (2.5) does not hold, that is, $p = P(G_0(\xi_p) = 1)$. Then by (2.7), $EG_0(\xi_p)I(0 < G_0(\xi_p) < 1) = 0$, which implies that $P(G_0(\xi_p) \in (0,1)) = 0$. Consequently,

$$P(G_0(\xi_p) = 0) = 1 - [P(G_0(\xi_p) \in (0,1)) + P(G_0(\xi_p) = 1)] = 1 - p.$$

But by the monotonicity of $G_0$,

$$P(A) = P(A \cap \{G_0(\xi_p) = 0\}) + P(A \cap \{G_0(\xi_p) = 1\})$$
$$\leq P(\{G_0((a - \varepsilon, a]) > \varepsilon\} \cap \{G_0(\xi_p) = 0\})$$
$$\quad + P(\{G_0(b, b + \varepsilon]) > \varepsilon\} \cap \{G_0(\xi_p) = 1\})$$
$$= P(\varnothing) + P(\varnothing) = 0,$$

which contradicts the fact that $P(A) > 0$. Hence, the claim is proved.

EXAMPLE 2.2. Let $\{Y_i\}_{i \in \mathbb{Z}}$ be a stationary homogeneous Markov process with transition probability function $P(\cdot; \cdot)$ and stationary distribution $\nu$. Let $X_i = H(Y_i), i \in \mathbb{N}$, where $H$ is a Borel measurable function. Suppose that

(2.8) $$|P(x; A) - P(y; A)| < 1$$

for all $x, y \in \mathbb{R}$ and $A \in \mathcal{B}(\mathbb{R})$, the Borel $\sigma$-field on $\mathbb{R}$. Then by (iii) on page 219 of Götze and Hipp (1983), conditions (C.2)–(C.4) hold with $\mathcal{D}_j = \sigma\langle Y_j \rangle$ and $X_{j,m}^\dagger = X_j$ for all $m \in \mathbb{N}, j \in \mathbb{Z}$ where (2.1)–(2.3) are satisfied with arbitrarily large positive real numbers $\alpha_0, \beta_0, \gamma_0$. For condition (C.5), we take $\mathcal{C}_j = \sigma\langle\{Y_i : i \neq j\}\rangle, j \in \mathbb{Z}$. Next, suppose that there exists a $\sigma$-finite measure $\mu$ such that $P(x; \cdot) \ll \mu$ for all $x \in \mathbb{R}$ and $\nu \ll \mu$. Write $f_0(x)$ and $f_1(x, y)$, respectively, for the density of $\nu$ and $P(x; \cdot)$ with respect to $\mu$. Also, let $A = H^{-1}((-\infty, \xi_p])$. Suppose that there exist sets $A_1 \subset A$ and $A_2 \subset A^c$ with $\mu(A_i) > 0$ for $i = 1, 2$ such that

(2.9) $\quad f_0(x) > 0, \qquad f_1(x, y) > 0 \qquad$ for all $x, y \in A_1 \cup A_2$.

We claim that condition (C.5) holds under (2.9).

To prove the claim, first note that the conditional distribution of $Y_0$ given $\{Y_i : i \neq 0\}$ is given by

$$(2.10) \quad P(Y_0 \in B | Y_i : i \neq 0) = \int_B g_0(Y_{-1}, Y_1; y) \mu(dy), \qquad B \in \mathcal{B}(\mathbb{R}),$$

where

$$g_0(y_{-1}, y_1; y_0)$$
$$= \begin{cases} \dfrac{f_1(y_{-1}, y_0) f_1(y_0, y_1)}{\int_{\mathbb{R}} f_1(y_{-1}, y) f_1(y, y_1) \mu(dy)}, & \text{if the denominator is positive,} \\ 0, & \text{otherwise.} \end{cases}$$

Relation (2.10) can be easily established by verifying the integral equation

$$\int_{D_k} P(Y_0 \in B | Y_i : i \neq 0) \, dP = P(D_k \cap \{Y_0 \in B\})$$

for all $B \in \mathcal{B}(\mathbb{R})$ and for all sets $D_k$ of the form $D_k = \bigcap_{0 < |i| \leq k} \{Y_i \in B_i\}$ for $k \in \mathbb{N}$, where $B_i \in \mathcal{B}(\mathbb{R})$ for all $i$.

Next note that by (2.9) and (2.10), for all $(y_{-1}, y_1) \in A_1 \times A_2$,

$$P(Y_0 \in A | y_{-1}, y_1) \geq P(Y_0 \in A_1 | y_{-1}, y_1) > 0,$$
$$P(Y_0 \in A^c | y_{-1}, y_1) \geq P(Y_0 \in A_2 | y_{-1}, y_1) > 0.$$

Consequently, $G_0(\xi_p) \equiv P(X_0 \leq \xi_p | Y_i, i \neq 0) = P(Y_0 \in A | Y_i, i \neq 0) = P(Y_0 \in A | Y_{-1}, Y_1) \in (0, 1)$ for all $(Y_{-1}, Y_1) \in A_1 \times A_2$. Since by (2.9), $P(Y_{-1} \in A_1, Y_1 \in A_2) > 0$, by the identity given in (2.7), $P(G_0(\xi_p) = 1) < p$. Hence, condition (C.5) follows.

EXAMPLE 2.3. Let $\{Y_j\}_{j \in \mathbb{Z}}$ be a stationary zero-mean unit variance Gaussian process with spectral density $f(\lambda)$ and let $X_j = H_0(Y_j), j \in \mathbb{Z}$ for some Borel measurable function $H_0 : \mathbb{R} \to \mathbb{R}$. For this example, we set $\mathcal{D}_j = \sigma\langle Y_j \rangle$ and $X_{j,m}^\dagger = X_j$ for all $j \in \mathbb{Z}$, $m \geq 1$. Then it is clear that condition (C.3) holds with an arbitrarily large $\beta_0 \in (0, \infty)$. Next, note that by Theorem V.6.8 of Ibragimov and Rozanov (1978), the strong mixing coefficient $\alpha(\cdot)$ of the Gaussian process $\{Y_i\}_{i \in \mathbb{Z}}$ satisfies $\alpha(n) = O(n^{-12-\delta})$ as $n \to \infty$ for some $\delta \in (0, 1)$ if and only if $f(\lambda)$ is of the form

$$(2.11) \qquad f(\lambda) = |p(\exp(\iota\lambda))|^2 w(\lambda), \qquad \lambda \in (-\pi, \pi]$$

where $p(z)$ is a polynomial with zeros on the unit circle $\{|z| = 1\}$ and where $w(\lambda)$ is a function that is bounded away from zero and is 12-times differentiable such that the 12th derivative satisfies a Hölder's condition of order $\delta$. Thus, condition (C.2) holds under (2.11). Also, note that by the arguments on page 220 of Götze and Hipp (1983), condition (C.4) holds with a $\gamma_0 > 12$.

To verify condition (C.5), we take $\mathcal{C}_j = \sigma\langle\{Y_i : i \neq j\}\rangle$. The arguments in Example 2.1 imply that condition (C.5) is not true if and only if

$$(2.12) \qquad P(G_0(\xi_p) = 1) = p \quad \text{and} \quad P(G_0(\xi_p) = 0) = 1 - p,$$

where $G_0(\xi_p) = P(X_0 \leq \xi_p | Y_i : i \neq 0) = P(H_0(Y_0) \leq \xi_p | Y_i : i \neq 0)$. Since the conditional distribution function of $Y_0$ given $(Y_i : i \neq 0)$ is normal and the sets $H_0^{-1}(-\infty, \xi_p]$ and $H_0^{-1}(\xi_p, \infty)$ both have positive probabilities under N(0, 1) [as $P(H_0(Y_0) \leq \xi_p) = p \in (0, 1)$], by the absolute continuity of normal distributions, $G_0(\xi_p) \in (0, 1)$ with probability one. Hence, (2.12) fails and, therefore, condition (C.5) holds. Thus, for the process $\{X_i\}_{i \in \mathbb{Z}}$ of this example, the Berry–Esseen theorem for the sample quantiles holds solely under condition (C.1) and (2.11). Note that (2.11) requires only a polynomial decay of the autocovariance function of the Gaussian process $\{Y_i\}_{i \in \mathbb{Z}}$.

2.3. *The Theorem.* We now state the main result of the paper.

THEOREM. *Suppose that conditions* (C.1)–(C.5) *hold. Then there exists a constant* $C \in (0, \infty)$ *such that for all* $n \geq 1$,

$$\sup_{x \in \mathbb{R}} |P(\sqrt{n}(\hat{\xi}_n - \xi_p) \leq x) - \Phi(x)| \leq \frac{C}{\sqrt{n}}$$

*where* $\hat{\xi}_n = F_n^{-1}(p)$ *and* $\xi_p = F^{-1}(p)$.

Thus, under conditions (C.1)–(C.5), the rate of normal approximation to the distribution of normalized sample quantiles of strongly mixing random variables is $O(n^{-1/2})$. This rate agrees with the standard rate available in the case of normalized sample mean and sample quantiles of independent random variables. As mentioned earlier, the $O(n^{-1/2})$ bound under dependence is rather surprising, as a similar bound in the case of the sample mean of strongly mixing random variables still remains elusive, even under an exponential decay of the mixing coefficient $\alpha(n)$ of (2.1).

The above theorem also extends the result of Reiss (1974) on the rate of normal approximation to the distribution of the sample quantiles of i.i.d. random variables, by allowing the random variables to be approximately strongly mixing. The method of proof employed here is very different from Reiss (1974) proof which heavily exploits the formula for the probability density functions of the sample quantiles of i.i.d. random variables. The same approach does not extend easily to the dependent case considered here as a similar formula for the density is not available for the general class of mixing processes. In contrast, our proof makes use of the characteristic function techniques of Götze and Hipp (1983) and Lahiri (1993, 1996), and some uniform bounds on the behavior of characteristic functions of indicator variables in the neighborhood of the population quantile $F^{-1}(p)$, which may be of some independent interest. See Lemma 3.3 in Section 3.

**3. Proofs.** In the proofs below, we write $C, C(\cdot)$ to denote generic constants with values in $(0, \infty)$ that may depend on the arguments (if any), but not on the variables, $n, x, y$. Also, unless otherwise mentioned, we take limits by letting $n \to \infty$. Let $\iota = \sqrt{-1}$. For any two real numbers $x, y$, let $x \wedge y = \min\{x, y\}$ and $x \vee y = \max\{x, y\}$.

By the definition of the sample quantile, for any $y \in \mathbb{R}$,

$$(3.1) \qquad P(F_n(y) > p) \leq P(\hat{\xi}_n \leq y) \leq P(F_n(y) \geq p).$$

Hence, we consider the sums $\sum_{i=1}^n I(X_i \leq y)$ for $y$ in a neighborhood of $\xi_p$ and study the rate of convergence of the upper and the lower bounds in (3.1). The first result gives an expansion for the log-characteristic function of a scaled sum of a transformed sequence $\{f_n(X_j)\}_{j \in \mathbb{Z}}$ of random variables in a neighborhood of the origin.

LEMMA 3.1. *For each $n \in \mathbb{N}$, let $f_n : \mathbb{R} \to [-1, 1]$ be a Borel measurable function such that $E f_n(X_i) = 0$, $E|f_n(X_i) - f_n(X_{i,k}^\dagger)| \leq Ck^{-\beta_0}$ for all $i \in \mathbb{Z}$, $k \in \mathbb{N}$ and*

$$(3.2) \qquad n^{-1} \operatorname{Var}\left(\sum_{i=1}^n f_n(X_i)\right) = 1.$$

*Let $W_{ni} \equiv f_n(X_i), i \in \mathbb{Z}, n \geq 1, S_n = n^{-1/2} \sum_{i=1}^n W_{ni}, H_n(t) = E \exp(\iota t S_n), t \in \mathbb{R}$. Also, let $\chi_{r,n}$ denote the $r$th cumulant of $S_n$. Then for any $\varepsilon \in (0, 1/4)$,*

$$\sup_{t \in A_n} \left| \log E \exp(\iota t S_n) - \sum_{r=2}^5 \frac{(\iota t)^r}{r!} \chi_{r,n} \right|$$

$$\leq C(\varepsilon) \left( \sup_{t \in A_n} |H_n(t)|^{-6} \right) \cdot n^{2\varepsilon(\alpha_0 \vee \beta_0)} \cdot \{n^{-1/2 - \alpha_0/4} + n^{-1/2 - \beta_0/4}\}$$

$$+ C(\varepsilon) n^{-1/2 - 6\varepsilon} \left( 1 + \sup_{t \in A_n} |\theta_{1n}(t)|^6 \right)$$

*for all $n \geq 1$, where $A_n = \{t \in \mathbb{R} : |t| \leq (\log n)^{1/2} (\log \log(n+1))^{1/4}\}$ and where $\theta_{1n}(t)$ is as defined in (3.7) below.*

PROOF. For any random variables $V_1, \ldots, V_p$, $p \in \mathbb{N}$, set

$$(3.3) \quad \begin{aligned} & \mathcal{K}_t(V_1, \ldots, V_p) \\ & = \frac{\partial}{\partial x_1} \cdots \frac{\partial}{\partial x_p} \log E \exp(\iota t S_n + x_1 V_1 + \cdots + x_p V_p)|_{x_1 = \cdots = x_p = 0}. \end{aligned}$$

Then using Taylor's expansion of the cumulant generating function "$\log E \exp(\iota t S_n)$" around $t = 0$, we get

$$(3.4) \qquad \left| \log E \exp(\iota t S_n) - \sum_{r=2}^5 \frac{(\iota t)^r}{r!} \chi_{r,n} \right| \leq \sum_{k=0}^{n-1} \sum^{(k)} |\mathcal{K}_{\eta t}(V_{j_1}, \ldots, V_{j_6})|$$

for any $t \in \mathbb{R}$ with $|E\exp(\iota t S_n)| > 0$, where $\eta \equiv \eta(t) \in [0,1]$, $V_j = tW_{nj}/\sqrt{n}$ and for a given $k$, the summation $\sum^{(k)}$ extends over $j_1, \ldots, j_6$ with maximal gap $k$. Note that by Lemma 3.1 of Lahiri (1996) (with $c_n = 1, t = 0$), for any $a_1, \ldots, a_r \in \mathbb{R}$, with $|a_j| \leq 1, r \geq 2$,

$$|\mathcal{K}_0(a_1 S_n, \ldots, a_r S_n)| \leq C(r) n^{-(r-2)/2} \sum_{k=0}^{n-1} k^{r-1}[\alpha(k/3) + \beta_*(k/3)]$$

(3.5)
$$\leq C(r) n^{-(r-2)/2},$$

provided $\alpha_0 > r, \beta_0 > r$, where $\beta_*(k) \equiv k^{-\beta_0}, k \in \mathbb{N}$.

Next, fix $\varepsilon \in (0, 1/4)$ and let $a_n = n^{1/4-\varepsilon}$. Then by (3.4) above and by Lemma 3.2 of Lahiri (1996) (with $c_n = 1$), as in the proof of his Lemma 3.6 [cf. (3.9), *op. cit.*],

$$\sum_{k=0}^{a_n} \sum^{(k)} |\mathcal{K}_{\eta t}(V_{j_1}, \ldots, V_{j_6})|$$

$$\leq \sum_{k=0}^{a_n} n(k+1)^5 C n^{-3}(1+|t|^6)$$

(3.6)
$$\times \{(1 + \theta_{1n}(\eta t))^6 + (1 + \theta_{2n}(\eta t))^6\}$$

$$\leq C(\varepsilon) a_n^6 n^{-2}(1+|t|^6)$$

$$\times [1 + |\theta_{1n}(t)|^6 + |H_n(t)|^{-6}\{n^{-3\alpha_0/4} + n \cdot n^{-3\beta_0/4}\}],$$

for all $t \in A_n$, where

(3.7) $\theta_{1n}(t) = \dfrac{1}{|H_n(t)|} \max\{|E\exp(S_I^{(l)})| : 1 \leq l \leq L, |I| \leq 4, I \subset \{1, \ldots, n\}\}$,

and where $\theta_{2n}(t) = |H_n(t)|^{-1}[L2^L\{\alpha(m) + n\beta_*(m)\} + \{\zeta_t(m)\}^L]$, $\zeta_t(k) = C|t| \times (n^{-1}k)^{1/2}$ [correcting for a typographical error in Lahiri (1996)] for $k \in \mathbb{N}$, $m = n^{3/4+\varepsilon}$ and $L = \log\log n$. Here, for $I \subset \{1, \ldots, n\}, l \geq 0$, $S_I^{(l)} \equiv \iota n^{-1/2} t \sum^{*(l)} W_{nj}$, where the summation $\sum^{*(l)}$ ranges over all $j \in \{1, \ldots, n\}$ such that $|j - i| > lm$ for all $i \in I$. Next using Lemma 3.3 of Lahiri (1996) with $K = L, m = 3Kn^{-\varepsilon}$ and $c_n = 1$, for each $k \in (a_n, n)$, as in the proof of (3.10), page 217 of Lahiri (1996) [correcting for the typo, where $(1 + \|t\|^{r/2})$ is replaced with $(1 + \|t\|^r)n^{-r/2}$], we get

$$\sum_{k=a_n+1}^{n-1} \sum^{(k)} |\mathcal{K}_{\eta t}(V_{j_1}, \ldots, V_{j_6})|$$

$$\leq C \dfrac{(1+|t|^6)}{n^3 |H_n(\eta t)|^6} L2^L$$

(3.8)
$$\times \left\{ \sum_{k=a_n+1}^{n} n(k+1)^5 [\alpha(kn^{-\varepsilon}) + n\beta_*(kn^{-\varepsilon}) + \zeta_t(3kn^{-\varepsilon})^L] \right\}$$
$$\leq C(\varepsilon)(1+|t|^6)|H_n(\eta t)|^{-6} \cdot n^{-2}$$
$$\times L2^L \{n^{\varepsilon \alpha_0} a_n^{(6-\alpha_0)} + n \cdot n^{\varepsilon \beta_0} \cdot a_n^{(6-\beta_0)} + n^{-\varepsilon L/4}\}$$
$$\leq C(\varepsilon) \cdot |H_n(\eta t)|^{-6} \cdot n^{2\varepsilon(\alpha_0 \vee \beta_0)} \cdot n^{-1/2} \cdot \{n^{-\alpha_0/4} + n \cdot n^{-\beta_0/4}\}$$

for all $t \in A_n$. Hence, the lemma follows from (3.4) and (3.6)–(3.8). □

REMARK. It is possible to obtain a bound on the difference between $E \exp(\iota t S_n)$ and its $s$th order Taylor expansion $\sum_{r=2}^{s} \frac{(\iota t)^r}{r!} \chi_{r,n}$ for an integer $s \geq 3$ by suitably modifying the arguments in the proof of Lemma 3.1. It can be shown that for a small $\delta > 0$, a bound of the order $O(n^{-1/2-\delta})$ on the difference is assured, if the strong-mixing exponent $\alpha_0$ satisfies

(3.9)
$$\alpha_0 > s + 4 + 9/(s-2).$$

By minimizing the right-hand side of (3.9), we get $s = 5$, which explains the reason behind considering the 5th order Taylor's expansion in Lemma 3.1.

LEMMA 3.2. *Let $W_{nj}$'s and $S_n$ be as in Lemma 3.1.*

(i) *Then for any $a \in (0, 1/2)$, there exists a constant $C_0 = C_0(\alpha_0, \beta_0, \gamma_0, a)$ such that for all $n \geq C_0$,*

$$|H_n(t)| \leq C_0 \left[ \exp\left( \frac{-t^2}{2} \left[ 1 - \frac{C_0}{(\log n)^2} \right] \right) + n^{1-a} \{n^{-a\alpha_0} + n^{-a\beta_0} + n^{-a\gamma_0}\} (\log n)^{C_0} \right]$$

*uniformly in $|t| \leq n^{(1-a)/2} (\log n)$.*

(ii) *There exist $\varepsilon_0 \in (0,1)$ and $C_1 = C_1(\alpha_0, \beta_0, \gamma_0) \in (0, \infty)$ such that for all $n \geq C_1$,*

$$|H_n(t)| \geq \varepsilon_0 \exp(-t^2/2) - C_1 \cdot (\log n)^{C_1} [n^{-\alpha_0/2} + n^{1/4} \cdot n^{-\beta_0/2} + n^{-\gamma_0/2}]$$

*for all $|t| \leq \varepsilon_0 \log n$.*

PROOF. Let $m = m_1(\log n)^{-2}$ and $m_1 = n^a(\log n)^{-6}$. Let $l, j_1, \ldots, j_l$ be the integers defined on page 218 of Lahiri (1996) with $I = \{1, \ldots, n\}$ and $I_1 = \{m_1+1, \ldots, n-m_1\}$. Also, let $\Gamma_k = \prod \{\exp(\iota t W_{nj}/\sqrt{n}) : j \in I, |j-j_k| \leq m_1\}$, $k = 1, \ldots, l$, and $B = \prod \{\exp(\iota t W_{nj}/\sqrt{n}) : j \in I, |j-j_k| > m_1 \text{ for all } k = $

$1, \ldots, l\}$. Then using the arguments leading to (3.11) of Lahiri (1996) (with $c_n = 1, R = 1$), for all $n \geq 1, t \in \mathbb{R}$, we get

$$|H_n(t)| = \left| E\left(\prod_{k=1}^{l} \Gamma_k\right) B \right|$$

(3.10)
$$\leq C \prod_{k=1}^{l} E|E(\Gamma_k | \mathcal{D}_j : j \neq j_k)|$$

$$+ C[l\alpha(m) + l\gamma(m) + \beta_*(m)\{m + n^{1/2}|t|\}].$$

Next, note that by (3.2) and the stationarity of $X_i$'s,

$$\left| m_1^{-1} \operatorname{Var}\left(\sum_{j=1}^{m_1} W_{nj}\right) - 1 \right|$$

$$\leq 2 \sum_{j=m_1+1}^{n-1} |EW_{n_1} W_{n(j+1)}| + \frac{4}{m_1} \sum_{j=1}^{n} j |EW_{n1} W_{n(j+1)}|$$

(3.11)
$$\leq C \left[ \sum_{j=m_1+1}^{n} \{\alpha(j/3) + \beta_*(j/3)\} + \frac{1}{m_1} \sum_{j=1}^{\infty} j \{\alpha(j/3) + \beta_*(j/3)\} \right]$$

$$\leq C(\alpha_0, \beta_0)[m_1^{-\alpha_0+1} + m_1^{-\beta_0+1} + m_1^{-1}].$$

Hence, by (3.12) and the arguments following it on page 219 of Lahiri (1996), and by (3.10) and (3.11) above, it follows that for all $n \geq C(\alpha_0, \beta_0)$,

$$\prod_{k=1}^{l} |E(\Gamma_k | \mathcal{D}_j : j \neq j_k)|$$

(3.12) $\leq C \exp\left(-\frac{t^2}{2}\{n^{-1} 2m_1 l(1 - C(\alpha_0, \beta_0) m_1^{-1}) - C n^{-3/2} l |t| m_1^{3/2}\}\right)$

$$\leq C \exp\left(-\frac{t^2}{2}[1 - C(a, \alpha_0, \beta_0)(\log n)^{-2}]\right)$$

for all $|t| \leq n^{(1-a)/2}(\log n)$, where in the second inequality, we have made use of the fact

(3.13) $\quad 2m_1 l = n[1 - O(n^{-1} m_1 + m_1^{-1} m)] \quad \text{as } n \to \infty.$

Hence, part (i) of the lemma follows from (3.10) and (3.12). Part (ii) can be proved by retracing the arguments on pages 221–222 of Lahiri (1996), with $c_n = 1$. We omit the routine details. $\square$

For the next lemma, let $\mathcal{C}$ be a sub-$\sigma$-algebra of $\mathcal{F}$, $G(y; \cdot) = P(X_1 \leq y | \mathcal{C})$, $A_1(y) = \{\omega : G(y; \omega) = 1\}$ and $A_2(y) = \{\omega : 0 < G(y; \omega) < 1\}, y \in \mathbb{R}$. Also, let

$g(y) = P(\{\omega : G(y; \omega) = 1\}) = P(A_1(y)), y \in \mathbb{R}$. Let $\Psi_a(t) = (ae^{\iota t} + 1 - a), t \in \mathbb{R}$ denote the characteristic function of a random variable $Y$ with $P(Y = 0) = 1 - a, P(Y = 1) = a, a \in (0, 1)$.

LEMMA 3.3. *If $g(\xi_p) < p$, then there exist $\delta, \varepsilon \in (0, 1)$ such that for all $t \in \mathbb{R}$,*

$$\sup_{|y - \xi_p| \leq \delta} E|E\{\exp(\iota t I(X_1 \leq y))|\mathcal{C}\}| \leq 1 - (1 - |\Psi_\varepsilon(t)|)\delta.$$

PROOF. By definition, for all $y \in \mathbb{R}$,

(3.14)
$$F(y) = P(X_1 \leq y) = E\{P(X_1 \leq y|\mathcal{C})\}$$
$$= \int_{A_1(y) \cup A_2(y)} G(y; \cdot) \, dP = \int_{A_2(y)} G(y; \cdot) \, dP + g(y).$$

Note that $A_1(y_1) \subset A_1(y_2)$ for all $y_1 < y_2$ and that $G(\cdot; \omega)$ is a valid distribution function for each $\omega \in \Omega$. We claim that:

(i) $g(\cdot)$ is nondecreasing, and
(ii) $g(\cdot)$ is right continuous on $\mathbb{R}$.

The first assertion is immediate. To prove (ii), note that for any sequence $y_n \downarrow y \in \mathbb{R}$,

$$\{\omega : G(y; \omega) = 1\} \subset \bigcap_{n \geq 1} \{\omega : G(y_n; \omega) = 1\}$$

[by the monotonicity of $G(y; \cdot)$ in $y$] while the reverse inclusion follows from the right continuity of $G(y; \omega)$ in $y$ for each $\omega$. Since $F(\xi_p) = p$, by (3.14),

(3.15)  $$g(\xi_p) \leq p.$$

Now suppose that $g(\xi_p) < p$. Then by the right continuity of $g(\cdot)$, there exists a $\delta_0 > 0$ such that

(3.16)  $$g(\xi_p + \delta_0) < p - \delta_0.$$

By (i), this implies that $g(y) < p - \delta_0$ for all $y < \xi_p + \delta_0$. Since $F$ is continuous at $\xi_p$, there exists a $0 < \delta_1 \leq \delta_0$ such that

(3.17)  $$F(\xi_p - \delta_1) > p - [\delta_0/2].$$

Next, write $A_3(y; \varepsilon) = \{\omega : \varepsilon < G(y; \omega) < 1 - \varepsilon\}, A_4 = \{\omega : G(y; \omega) \geq 1 - \varepsilon\}, \varepsilon \in (0, 1), y \in \mathbb{R}$. Note that $A_4(y; \varepsilon) \downarrow A_1(y)$ as $\varepsilon \downarrow 0$ and $A_4(y; \varepsilon) \subset A_4(y + h; \varepsilon)$ for all $y \in \mathbb{R}, h > 0, \varepsilon \in (0, 1)$. In particular, for any $y \in \mathbb{R}$,

(3.18) $$\lim_{\varepsilon \downarrow 0} \int_{A_4(y; \varepsilon)} G(y; \cdot) \, dP = \int_{A_1(y)} G(y; \cdot) \, dP = P(A_1(y)) = g(y).$$

Hence, by (3.14) and (3.16)–(3.18), there exists $0 < \varepsilon < \delta_0/8$ such that for all $y \in (\xi_p - \delta_1, \xi_p + \delta_1)$,

$$
\begin{aligned}
P(A_3(y;\varepsilon)) &\geq \int_{A_3(y;\varepsilon)} G(y;\cdot)\,dP \\
&= \int G(y;\cdot)\,dP - \int G(y;\cdot)I(G(y;\cdot) \leq \varepsilon)\,dP \\
&\quad - \int G(y;\cdot)I(G(y;\cdot) \geq 1-\varepsilon)\,dP \\
&\geq F(y) - \varepsilon - \int G(\xi_p + \delta_1;\cdot)I(G(\xi_p + \delta_1;\cdot) \geq 1-\varepsilon)\,dP \\
&\geq F(\xi_p - \delta_1) - \varepsilon - [g(\xi_p + \delta_1) + \varepsilon] \\
&\geq [p - \delta_0/2] - 2\varepsilon - [p - \delta_0] \\
&= \delta_0/4.
\end{aligned}
$$
(3.19)

Next, writing $\Psi_\varepsilon(t) = |\varepsilon e^{\iota t} + 1 - \varepsilon|, t \in \mathbb{R}$, and $G(y) = G(y;\cdot)$ (for notational simplicity), by (3.19), for $y \in (\xi_p - \delta_1, \xi_p + \delta_1)$, we have

$$
\begin{aligned}
E|E(\exp(\iota t I(X_1 \leq y))|\mathcal{C})| \\
&= E|G(y)e^{\iota t} + (1 - G(y))| \\
&= E|1 - 4G(y)(1 - G(y))\sin^2(t/2)|^{1/2} \\
&\leq P(A_3^c(y;\varepsilon)) + EI_{A_3(y;\varepsilon)} \cdot |1 - 4G(y)(1 - G(y))\sin^2(t/2)|^{1/2} \\
&\leq P(A_3^c(y;\varepsilon)) + |1 - \varepsilon(1-\varepsilon)\sin^2(t/2)|^{1/2} P(A_3(y;\varepsilon)) \\
&= 1 - (1 - |\Psi_\varepsilon(t)|)P(A_3(y;\varepsilon)) \\
&\leq 1 - (1 - |\Psi_\varepsilon(t)|)\delta_0/4. \qquad \square
\end{aligned}
$$
(3.20)

PROOF OF THE THEOREM. First, we shall show that

$$
\Delta_n^* \equiv \sup_{|x| \leq \log n} |P(\sqrt{n}(\hat{\xi}_n - \xi_p) \leq x) - \Phi(x/\tau_\infty(p))|
$$
(3.21)
$$
= O(n^{-1/2}).
$$

To prove this, we apply inequality (3.1) with $y = x_n$, where for $x \in \mathbb{R}$, we set $x_n = \xi_p + n^{-1/2}x$. Let $\sigma_n^2(x) = n\operatorname{Var}(F_n(x))$ and $S_n(x) = \sqrt{n}(F_n(x) - F(x))/\sigma_n(x)$, $x \in \mathbb{R}$. By the smoothing inequality [cf. Lemma 2, page 538 of Feller (1971)],

$$
\Delta_n \equiv \sup_{|x| \leq \log n} \left| P(F_n(x_n) \leq p) - \Phi\left(\frac{\sqrt{n}(p - F(x_n))}{\sigma_n(x)}\right) \right|
$$

$$
\begin{aligned}
(3.22) \quad &\leq \sup_{|x|\leq \log n} \sup_{y\in\mathbb{R}} |P(S_n(x_n) \leq y) - \Phi(y)| \\
&\leq \sup_{|x|\leq n^\delta \log n} \left[ \frac{1}{\pi} \int_{-\kappa\sqrt{n}}^{\kappa\sqrt{n}} |E\exp(\iota t S_n(x_n)) - e^{-t^2/2}||t|^{-1}\,dt + \frac{C}{\kappa\sqrt{n}} \right],
\end{aligned}
$$

where $\kappa \in (0,\infty)$ is a constant (independent of $x$), to be specified later and where $C \in (0,\infty)$ is a constant independent of $n$, $x$.

Next note that

$$\sup\left\{ \sum_{j=n}^{\infty} |\mathrm{Cov}(I(X_1 \leq x), I(X_{j+1} \leq x))| : x \in \mathbb{R} \right\}$$

$$\leq C \sum_{j=n}^{\infty} [\alpha(j/3) + \beta(j/3))$$

$$\to 0 \quad \text{as } n \to \infty.$$

Since $\sigma_\infty^2(\xi_p) > 0$ and $F$ is continuous at $\xi_p$, by the above fact, there exists $\delta_* \in (0, \delta_0)$ such that

$$(3.23) \qquad \liminf_{n\to\infty} \inf\{\sigma_n^2(\xi_p + x) : |x| \leq \delta_*\} > \sigma_\infty^2(\xi_p)/2,$$

where $\delta_0$ is as in Lemma 3.3. Let $\mathcal{N}_p = \{x : |x - \xi_p| \leq \delta_*\}$, and let $\chi_{r,n}(x)$ denote the $r$th cumulant of $S_n(x) \equiv \sqrt{n}(F_n(x) - F(x))/\sigma_n(x), x \in \mathcal{N}_p$. Then it is easy to check that for any $i \in \mathbb{Z}$, any $x \in \mathbb{R}$ and any $\varepsilon_0 > 0$,

$$
\begin{aligned}
& E|I(X_i \leq x) - I(X_{i,k}^\dagger \leq x)| \\
(3.24) \quad & \leq P(x - \varepsilon < X_i < x + \varepsilon) + P(|X_i - X_{i,k}^\dagger| \geq \varepsilon) \\
& \leq [F(x+\varepsilon) - F(x-\varepsilon)] + \varepsilon^{-1} E|X_i - X_{i,k}^\dagger|.
\end{aligned}
$$

Hence, by conditions (C.1) and (C.3), there exists $C > 0$ and $\delta_{**} \in (0, \delta_*)$ such that for all $|x - \xi_p| < \delta_{**}$, and $i \in \mathbb{Z}, k \in \mathbb{N}$ [with $\varepsilon = k^{-\beta_0}$ in (3.24)],

$$(3.25) \qquad E|I(X_i \leq x) - I(X_{i,k}^\dagger \leq x)| \leq C k^{-\beta_0}.$$

For notational simplicity, without loss of generality, we shall set $\delta_{**} = \delta_*$. Also, let $W_{ni}(x) = [I(X_i \leq x) - F(x)]/\sigma_n(x), i \in \mathbb{Z}, n \geq 1, x \in \mathcal{N}_p$. Then by (3.23), (3.25) and condition (C.3), $\{W_{ni}(x) : i \in \mathbb{Z}\}_{n\geq 1}$ satisfies the conditions of Lemma 3.1 uniformly in $x \in \mathcal{N}_p$. By Lemmas 3.1 and 3.2(ii) above and the induction arguments used in the proof of Lemma 3.28 of Götze and Hipp (1983), it follows that there exists $\varepsilon_1 = \varepsilon_1(\alpha_0, \beta_0) \in (0, 1/4)$ such that for all $0 < \varepsilon \leq \varepsilon_1$,

$$\sup_{x\in\mathcal{N}_p} \sup_{t^2 \leq \log n} \left| \log E\exp(\iota t S_n(x)) - \sum_{r=2}^{5} \frac{(\iota t)^r}{r!} \chi_{r,n}(x) \right|$$

$$
\begin{aligned}
(3.26) \quad &\leq C(\varepsilon)[n^{-1/2-6\varepsilon} + \{\varepsilon_0 \exp(-(\sqrt{\log n})^2/2)\}^{-6} \\
&\qquad \times n^{2\varepsilon(\alpha_0 \vee \beta_0)}\{n^{-1/2-\alpha_0/4} + n^{-1/2-\beta_0/4}\}] \\
&\leq C(\varepsilon, \varepsilon_0) n^{-1/2-C(\varepsilon)}
\end{aligned}
$$

for all $n \geq C_1$, where $C_1$ is as in Lemma 3.2(ii). By arguments in the proofs of Lemma 3.33 of Götze and Hipp (1983) and of Lemma 9.7 of Bhattacharya and Range Rao (1976) and by (3.5) and (3.26), we have

$$
\sup_{x \in \mathcal{N}_p} \int_{t^2 \leq \log n} |E \exp(\iota t S_n(x)) - e^{-t^2/2}||t|^{-1}\, dt
$$

$$
\leq \sup_{x \in \mathcal{N}_p} \int_{t^2 \leq \log n} \left| E \exp(\iota t S_n(x)) \right.
$$

$$
\left. - e^{-t^2/2}\left(1 + \sum_{r=3}^{5}(\iota t)^r (r!)^{-1}\chi_{r,n}(x)\right)\right| |t|^{-1}\, dt
$$

$$
+ \sup_{x \in \mathcal{N}_p} \int_{t^2 \leq \log n} \left| e^{-t^2/2}\left(\sum_{r=3}^{5}(\iota t)^r (r!)^{-1}\chi_{r,n}(x)\right)\right| |t|^{-1}\, dt
$$

$$
\leq C n^{-1/2}
$$

for all $n \geq C_1$. Since

$$
\int_{t^2 > \log n} e^{-t^2/2}\, dt = o(n^{-1/2}),
$$

to show that $\delta_n$ of (3.22) is $O(n^{-1/2})$, it remains to show that

$$
(3.27) \quad \sup_{x \in \mathcal{N}_p} \int_{(\log n)^{1/2} < |t| < \kappa n^{1/2}} |E\exp(\iota t S_n(x))||t|^{-1}\, dt = O(n^{-1/2}).
$$

To this end, we split the set of $t$-values in (3.27) into the sets $B_{1n} = \{t \in \mathbb{R} : (\log n)^{1/2} \leq |t| \leq n^{7/16}\}$ and $B_{2n} = \{t \in \mathbb{R} : n^{7/16} < |t| < \kappa n^{1/2}\}$. Then using Lemma 3.2(i) with $W_{nj} = W_{nj}(x), x \in \mathcal{N}_p$ with $a = 1/8$, we have

$$
\sup_{x \in \mathcal{N}_p} \int_{B_{1n}} |E\exp(\iota t S_n(x))||t|^{-1}\, dt
$$

$$
\leq 2C_0 \int_{(\log n)^{1/2}}^{n^{7/16}} \exp(-t^2/2) \cdot \exp(C_0 t^2 (\log n)^{-2}/2)|t|^{-1}\, dt
$$

$$
+ 2C_0 \left(\int_{(\log n)^{1/2}}^{n^{7/16}} |t|^{-1}\, dt\right)(\log n)^{C_0} \cdot n^{-1/2-1/8}
$$

$$
(3.28) \quad \leq 2C_0 \left[\exp(C_0/2)\int_{(\log n)^{1/2}}^{\log n} \exp(-t^2/2)|t|^{-1}\, dt \right.
$$

$$+ \int_{\log n}^{n^{7/16}} \exp(-t^2/4)|t|^{-1}\, dt \Bigg]$$

$$+ 2C_0(\log n)^{C_0+1} n^{-1/2-1/8}$$

$$= o(n^{-1/2}).$$

Next, note that $(1-u)^{1/2} \leq 1 - u/2$ for all $0 < u < 1$. Hence, there exists $\kappa_0 = \kappa_0(\delta, \varepsilon) \in (0, \infty)$, depending on $\delta$ and $\varepsilon$ of Lemma 3.3 such that for $|t| \leq \kappa_0$,

$$\begin{aligned}
1 - (1 - |\Psi_\varepsilon(t)|)\delta &= (1-\delta) + \delta(1 - 4\varepsilon(1-\varepsilon)\sin^2(t/2))^{1/2} \\
&\leq (1-\delta) + \delta(1 - 2\varepsilon(1-\varepsilon)\sin^2(t/2)) \\
&\leq 1 - C(\varepsilon, \delta) t^2.
\end{aligned}$$

Also, note that for a bounded random variable $y$ and $\sigma$-fields $\mathcal{G}_1 \subset \mathcal{G}_2 \subset \mathcal{F}$, $E(Y|\mathcal{G}_1) = E\{E(Y|\mathcal{G}_2)|\mathcal{G}_1\}$ a.s. $(P)$. Hence, setting $\kappa = \kappa_0$ in $B_{2n}$, and using (3.10) (with $a = 1/8$), we have

$$\sup_{x \in \mathcal{N}_p} \int_{B_{2n}} |E \exp(\iota t S_n(x))| |t|^{-1}\, dt$$

$$\leq \sup_{x \in \mathcal{N}_p} \int_{B_{2n}} \prod_{k=1}^{l} E|E(\Gamma_k(x)|\mathcal{D}_j : j \neq j_k)| |t|^{-1}\, dt$$

$$+ C[n^{1-a}(n^{-a(\alpha_0 \wedge \beta_0 \wedge \gamma_0)}) \cdot (\log n)^{C(\alpha_0, \beta_0, \gamma_0)}]$$

$$\leq 2\log n \cdot \sup\left\{ \prod_{k=1}^{l} E|E(\Gamma_k(x)|\mathcal{C})| : t \in B_{2n}, x \in \mathcal{N}_p \right\} + o(n^{-1/2})$$

(3.29)
$$= 2\log n \cdot \sup\{E|E(\exp(\iota t I(X_1 \leq x)/\sqrt{n})|\mathcal{C})| : t \in B_{2n}, x \in \mathcal{N}_p\}^l$$

$$+ o(n^{-1/2})$$

$$\leq (2\log n) \cdot \sup_{t \in B_{2n}} \{1 - C(\varepsilon, \delta) \cdot t^2/n\}^l + o(n^{-1/2})$$

$$= O(\log n \cdot \exp(-C(\varepsilon, \delta) \cdot n^{-1/8} \cdot l)) + o(n^{-1/2})$$

$$= o(n^{-1/2}),$$

where $l = n/(2m_1)(1 + o(1)) = n^{1-a}(\log n)^6(1 + o(1))$ (with $a = 1/8$) and where the variables $\Gamma_k(x)$'s are defined as in the proof of Lemma 3.2 with $W_{jk} = W_{jk}(x), x \in \mathcal{N}_p$. Hence, (3.27) follows and by (3.1) and the absolute continuity of the limiting normal distribution, (3.21) follows.

Next, note that

$$\sup_{x \leq -\log n} |P(\sqrt{n}(\hat{\xi}_n - \xi_p) \leq x) - \Phi(x/\tau_\infty(p))|$$

$$
\begin{aligned}
(3.30) \quad & \leq P(\sqrt{n}(\hat{\xi}_n - \xi_p) \leq -\log n) + \Phi(-\log n/\tau_\infty(p)) \\
& \leq \Delta_n^* + 2\Phi(-\log n/\tau_\infty(p)) \\
& = O(n^{-1/2})
\end{aligned}
$$

and similarly,

$$
(3.31) \quad \sup_{x \geq \log n} |P(\sqrt{n}(\hat{\xi}_n - \xi_p) \leq x) - \Phi(x/\tau_\infty(p))| = O(n^{-1/2}).
$$

Hence, the theorem follows from (3.21), (3.30) and (3.31). $\square$

Department of Statistics  
Texas A&M University  
TAMU 3143  
College Station, Texas 77843-3143  
USA  
E-mail: snlahiri@tamu.edu

Department of Mathematics and Statistics  
Wright State University  
3640 Colonel Glenn HWY  
Dayton, Ohio 45435-0001  
USA  
E-mail: shuxia.sun@wright.edu